\documentclass{article}
\usepackage{amsfonts}
\usepackage{enumerate}
\usepackage{mathtools}

\newtheorem{co}{Corollary}
\newtheorem{lema}{Lemma}
\newtheorem{teo}{Theorem}

\newcommand{\iz}[1]{{\flushleft \bf #1}}
\title{Levinson Theorem for Differential Equations with Piecewise Constant Argument Generalized}
\author{Samuel Castillo\thanks{ Supported by DIUBB 074108 1/R}\\
Departamento de Matem\'atica.
Facultad de Ciencias.\\ Universidad del B\'{\i}o-B\'{\i}o. Casilla
5-C. Concepci\'on. Chile.\\ \tt{scastill@ubiobio.cl}
\and Wilfred Flores\\
Facultad de Ingenier\'{\i}a.
Universidad de Talca.\\ Campus Curic\'o. Camino Los Niches Km 1, Curic\'o. Chile.\\
\tt{wflores51@gmail.com}}
\date{}
\begin{document}
\maketitle
\begin{abstract}
In this work, it is presented an adaptation of an
asymptotic theorem of N. Levinson of 1948, to differential equation
with piecewise constant argument generalized, which were introduced
by M. Akhmet in 2007. The N. Levinson's theorem which is adapted is that dealt by M. S. P. Eastham in his work which is present in this bibliography. The more
relevant hypotheses of this theorem are highlighted an it is
established a version of this theorem with these hypotheses for
ordinary differential equations. Such a version is that which is
adapted to differential equation with piecewise constant argument
generalized. The adaptation is proved by mean the Banach fixed point
where contractive operator is built form a suitable version of the
constant variation formula.
\end{abstract}
\maketitle
\iz{2010 AMS Subject Class:} 34A38, 34C27, 34D09, 34D20.
\iz{Key words:} Differential Equations, Piecewise constant argument, Asymptotic Formula, Levinson Theorem.

\section{Introduction}
Let $\mathbb{N},\mathbb{N}_0,\mathbb{Z},\mathbb{R}, \mathbb{C}$ the sets of the positive integer numbers, of the nonnegative integer numbers, of the integer numbers, of the real numbers, of the complex numbers, respectively. Let $N \in \mathbb{N}$. We will consider $\mathbb{K}^N:={\cal M}_{N \times 1}(\mathbb{K})$, i.e, as column vectors, where $\mathbb{K} \in \{\mathbb{R},\mathbb{C}\}$. We will denote by $|\cdot|$ to the Euclidean norm for $\mathbb{K}^N$. ${\cal M}_{N}(\mathbb{C})$ will denote the $N \times N$ matrix with complex entries. On ${\cal M}_{N}(\mathbb{C})$,  $\|\cdot\|$ will denote the classic norm operator which is defined for $A \in {\cal M}_{N}(\mathbb{C})$, by $\displaystyle \|A\|=\sup_{v \in \mathbb{C}^N-\{0\}}\frac{|Av|}{|v|}.$

A Differential Equations with Piecewise Constant Argument Generalized (DEPCAG) is a differential equation of the form
\begin{equation}
\label{ieq}
\frac{dx}{dt}=f(t,x(t),x(\gamma (t))),
\end{equation}
where $\gamma:[t_0,+\infty[ \to [t_0,+\infty[$ is such that there is an strictly increasing secuence $\displaystyle (t_n)_{n=0}^{+\infty}$ with $t_n \to +\infty$ as $n \to +\infty$ and $\gamma([t_{n},t_{n+1}[)=\{\xi_n\}$ for all $n \in \mathbb{R}$. A function  $x=x(t)$ is understood to be  {\it solution of the DEPCAG (\ref{ieq})} if:
\begin{enumerate}
\item $x$ is continuous on $[t_0,+\infty[$;
\item the derivative $\displaystyle \frac{dx}{dt}$ of $x$ with the possible exception in $t=t_n$ for $n \in \mathbb{N}_0$, where is unilateral derivative exists;
\item $x$ is a solution of (\ref{ieq}) with the possible exception in $t=t_n$ for all $n \in \mathbb{N}_0$.
\end{enumerate}
Notice that (\ref{ieq}) is an ordinary differential equation in each interval $[t_n,t_{n+1}[$ for all $n \in \mathbb{N}_0$, but the leaps be between those intervals creates a difference system of the form
\[
\begin{array}{rcl}
x(\xi_n)&=&x(t_n)+\int_{t_{n}}^{\xi_n} f(\zeta,x(\zeta),x(\xi_n))d\zeta,\\
\\
x(t_{n+1})&=&x(\xi_n)+\int_{\xi_{n}}^{t_{n+1}} f(\zeta,x(\zeta),x(\xi_n))d\zeta,\\
\end{array}
\]
for all $n \in \mathbb{N}_0$.

The name generalized in those equations, is explained by the inclusion of the
differential equations with piecewise argument whose study seems to be started by K. Cooke and J. Wiener in \cite{CoWi84,Wi83,Wi93}.
They consider equations of the form $\ref{ieq}$ with
$\displaystyle \gamma (t)=\left[ t \right] $ or $\displaystyle \gamma
(t)=2\left[ {\frac{t+1}{2}} \right] $, where $[\cdot ]$ is the function assigns to each real number, the greater integer number less that it. The first known generalization was made by M. Akhmet \cite{A}.

In this work we establish a version of the asymptotic Levinson's
theorem  N. Levinson \cite[(1948)]{Lev} (see Eastham \cite{Ea}) for DEPCAG
\begin{equation}
\label{idepca2}
\frac{dy}{dt}=A(t)y(t)+B(t)y(\gamma(t))+F(t,z(g(t))),\;t\in
[t_0,+\infty[,
\end{equation}
where $A(t)$, $B(t)$  are matrices in ${\cal
M}_N(\mathbb{C})$, whose coefficients are locally integrable
functions of $t$, $F:[t_0,+\infty[ \times \mathbb{C}^N \to \mathbb{C}^N$ is a locally integrable function in the first variable such that there is
$\eta:[t_0,+\infty[ \to \mathbb{R}_0^{+}$ such that
\begin{equation}
\label{nl-lpchzt}
\begin{array}{rcl}
|F(t,\hat{a})-F(t,\hat{b})|&\leq&\eta(t)|\hat{a}-\hat{b}|\\
\\
F(t,0)&=&0,\\
\end{array}
\end{equation}
$\eta(t)$ will satisfy an $L^1$ type condition which will be given below(see (\ref{nl-l1})), $\gamma$ is a piecewise constant argument defined as above and $g:[t_0,+\infty[ \to [t_0,+\infty[$ satisfy 
\[
g([t_{n},t_{n+1}[) \subseteq [t_{n},t_{n+1}[,
\]
for all $n \in \mathbb{N}_0$.

Notice that $g$ could be an piecewise constant argument.

DEPCAG (\ref{idepca2}) will be seen as a
perturbation of the DEPCAG
\begin{equation}
\label{ihomo2} 
\frac{dz}{dt}=A(t)z+B(t)z(\gamma(t))
\end{equation}
which will have dichotomic condition similar to those given in the asymptotic Levinson theorem considered here.

An asymptotic result for DEPCAG is given by M. Akhmet \cite{A03}. He considers
the DEPCA,
\begin{equation}
\label{akh01}
\frac{dy}{dt}=C_0y+f(t,x(t),x(\gamma(t))),
\end{equation}
as a perturbation of the autonomous ordinary differential
equation
\begin{equation}
\label{akh02}
\frac{dx}{dt}=C_0x,
\end{equation}
where $\gamma$ is defined as above and the following hypotheses are
given
\begin{subequations}
\begin{equation}
\label{H1}
\exists L>0: \|f(t,x_1,y_1)-f(t,x_2,y_2)\| \leq
L(\|x_2-x_1\|+\|y_2-y_1\|)\;\mbox{y}\;f(t,0,0)=0;
\end{equation}
\begin{equation}
\label{H2}
\exists \overline{t}: 0 <t_{n+1}-t_n \leq \overline{t};
\end{equation}
\begin{equation}
\label{H3}
\begin{array}{l}
\exists M,m>0: m \leq \|e^{C(t-s)}\| \leq M,\forall t,s \in
[t_n,t_{n+1}]; ML\overline{t}e^{ML\overline{t}}<1;\\
2ML\overline{t}<1;
M^2L\overline{t}\left(\frac{ML\overline{t}e^{ML\overline{t}}+1}{1-ML\overline{t}e^{ML\overline{t}}}+ML\overline{t}e^{ML\overline{t}}\right)<m;
\end{array}
\end{equation}
\begin{equation}
\label{H4}
\exists \eta:\mathbb{R}^+\to [0,L]: \|f(t,x_1,y_1)-f(t,x_2,y_2)\| \leq \eta(t)(\|x_2-x_1\|+\|y_2-y_1\|);
\end{equation}
\begin{equation}
\label{H5}
\ell_0:=\int_0^{+\infty}
t^{m_{\beta}+m_{\alpha}-2}e^{t(\beta-\alpha)}\eta(t)dt<+\infty,
\end{equation}
where $\lambda_1,\ldots,\lambda_p$ are the characteristic values of $C_0$, $\displaystyle \alpha=\min_{j=1,\ldots,p}Re(\lambda_j)$, $\displaystyle \beta=\max_{j=1,\ldots,p}Re(\lambda_j)$, $m_{\alpha}$ and $m_{\beta}$  the maximum orders of the characteristic values of $C$ with real part equal to $\alpha$ and $\beta$ respectively, for all $t \in \mathbb{R}^{+}$.
\end{subequations}

Then, Akhmet \cite{A03} provides the following result.
\begin{teo}
\label{akh03} Assume that (\ref{H1})-(\ref{H5}) hold. Then,  for every solution
$y=y(t)$ of the DEPCAG (\ref{akh01}) has a representation
\begin{equation}
y(t)=e^{C_0t}[c+w(t)],
\end{equation}
where $w(t) \to 0$ as $t \to +\infty$ and $c \in \mathbb{R}^N$.
\end{teo}

That result, consider a DEPCAG as a perturbation of an autonomous
ordinary differential equation (\ref{akh02}). In Theorem
\ref{akh03}, it is seen the perturbation of the whole fundamental
matrix of (\ref{akh02}). In our result we see the perturbation of
only one dimension of the solution space of (\ref{ihomo2}), although
our non perturbed equation is already a DEPCAG. with the elements
for our proof  The last section is devoted to the main result and
its proof.

\section{Preliminaries}

In this section sets the definitions and facts to present the main
result.

First we find the Cauchy matrix for system (\ref{ihomo2}). 

We ask the condition
\begin{subequations}
\begin{equation}
\label{24a} D_n(t)=I+\int_{\xi_n}^{t} X(\xi_n,u)B(u)du\;\mbox{is
invertible}
\end{equation}
\end{subequations}
for all $t \in [t_n,t_{n+1}]$ and $n \in \mathbb{N}_0$, where
$X(t,s)=X(t)X(s)^{-1}$ and $X$ is a fundamental matrix for the system
\begin{equation}
\label{ihomo} \frac{dx}{dt}=A(t)x.
\end{equation}

Let
\[
\Phi(n)=H(n-1)H(n-2)\cdots H(0),
\]
where
\begin{equation}
H(n)=X(t_{n+1},\xi_n)D_n(t_{n+1})D_n(t_n)^{-1}X(\xi_n,t_n).
\end{equation}

for all $n \in \mathbb{N}_0$ and for all $t \in [t_n,t_{n+1}]$.

Given $s,t \in [t_{0},+\infty[$, we have that if $Z:[t_0,+\infty[
\to {\cal M}_N(\mathbb{C})$ is a matrix function such that
$Z(s,s)=I$ and for all $z_0 \in \mathbb{C}^N$, $z=Z(\cdot,s)z_0$ is
a solution of (\ref{ihomo2}) in the sense given in the introduction
for equation (\ref{ieq}) such that $z(s)=z_0$, then
\begin{equation}
\label{mfhomo2}
\begin{array}{rcl}
Z(t,s)&=&X(t,\gamma(t))D_{k_t}(t)D_{k_t}(t_{k_t})^{-1}X(\gamma(t),t_{k_t})\\
\\
&\times&\Phi(k_t)\Phi(k_s+1)^{-1}\\
\\
&\times&X(t_{k_s+1},\gamma(s))D_{k_s}(t_{k_s+1})D_{k_s}(s)^{-1}X(\gamma(s),s).\\
\\
\end{array}
\end{equation}
% Then,
%\[
%Z(t,s)=X(t,\xi_n)\left[I+\int_{\xi_n}^{t} X(\xi_n,u)B(u)du\right]\left[I+\int_{\xi_n}^{s} X(\xi_n,u)B(u)du\right]^{-1}X(\xi_n,s),
%\]

%{\tt For simplicity, we will consider} $\xi_n=t_n$, for all $n \in
%\mathbb{N}_0$. It will be seen that restriction does not imply an
%important loss of generality.

Assume that there are an unitary vector $\hat{e}$  and locally integrable functions $\hat{\lambda}:[t_0,+\infty[ \to \mathbb{C}$ and $\hat{\lambda}_d:[t_0,+\infty[ \to \mathbb{C}$ such that
\begin{subequations}
%\begin{align}
\begin{equation}
\label{25a}
\begin{array}{rcl}
Z(t,s)\hat{e}&=&\tilde{e}(t,s)\hat{e},\\
\end{array}
\end{equation}
{where}
\begin{equation}
\label{25b}
\begin{array}{rcl}
\tilde{e}(t,s)&=&e^{\int_s^t \hat{\lambda}(\xi)d\xi}\\
%&\times&\left[1+\int_{\gamma(s)}^{s}e^{-\int_{\gamma(s)}^{u}\hat{\lambda}(\xi)d\xi}\hat{\lambda}_{d}(u)du\right]^{-1}\\
\\
&\times&\left(1+\int_{\gamma(t)}^te^{\int_{\gamma(t)}^u\hat{\lambda}(\xi)d\xi}\hat{\lambda}_d(u)du\right)\left(1+\int_{\gamma(t)}^{t_{k_t}}e^{\int_{\gamma(t)}^{u}\hat{\lambda}(\xi)d\xi}\hat{\lambda}_d(u)du\right)^{-1}\\
\\
&\times&\left(\prod_{j=k_s+1}^{k_t}\left[1+\int_{\xi_j}^{t_{j+1}}e^{\int_{\xi_j}^{u}\hat{\lambda}(\xi)d\xi}\hat{\lambda}_d(u)du\right]
\left[1+\int_{\xi_j}^{t_{j}}e^{\int_{\xi_j}^{u}\hat{\lambda}(\xi)d\xi}\hat{\lambda}_d(u)du\right]^{-1}\right)\\
\\
&\times&\left(1+\int_{\gamma(s)}^{t_{k_s+1}}e^{\int_{\gamma(s)}^{u}\hat{\lambda}(\xi)d\xi}\hat{\lambda}_d(u)du\right)\left(1+\int_{\gamma(s)}^se^{\int_{\gamma(s)}^u\hat{\lambda}(\xi)d\xi}\hat{\lambda}_d(u)du\right)^{-1}\\
\\
\end{array}
\end{equation}
%\end{align}
\end{subequations}
and $k_t \in \mathbb{N}_0$ is defined for all $t \in [t_0,+\infty[$ such that $t \in [t_{k_t},t_{k_t+1}[$.

Condition (\ref{24a}) implies the invertibility of $Z(t,s)$ which
allows to define $Z(t,s)=Z(s,t)^{-1}$ for $t<s$.

Now, we use the variation of constants formula for the DEPCAG
\begin{equation}
\label{nh}
\frac{d\psi}{dt}=A(t)\psi(t)+B(t)\psi(\gamma(t))+f(t).
\end{equation}

It  can be written as,
\[
\psi(t)=Z(t,\gamma(t))\psi(\gamma(t))+\int_{\gamma(t)}^t X(t,s)f(s)ds.
\]
and
\[
\psi(\gamma(t))=Z(\gamma(t),t_{k_t})\psi(t_{k_t})+\int_{t_{k_t}}^{\gamma(t)}
X(\gamma(t),s)f(s)ds.
\]
Then,
\[
\psi(t)=Z(t,\gamma(t))\left[Z(\gamma(t),t_{k_t})\psi(t_{k_t})+\int_{t_{k_t}}^{\gamma(t)}
X(\gamma(t),s)f(s)ds\right]+\int_{\gamma(t)}^t X(t,s)f(s)ds,
\]
i.e.,
\[
\psi(t)=Z(t,t_{k_t})\psi(t_{k_t})+\int_{t_{k_t}}^{t}
\Gamma(t,s)f(s)ds.
\]
So,
\[
\psi(t_{n+1})=H(n)x(t_n)+\int_{t_n}^{t_{n+1}}
\Gamma(t_{n+1},s)f(s)ds,
\]
for all $n \in [n_0,+\infty [\cap \mathbb{Z}$,
where
\begin{equation}
\label{grpinto} \Gamma(t,s)= \left\{\begin{array}{rcl}
Z(t,\gamma(t))X(t,s)&\mbox{if}&t_{k_t} \leq s \leq \gamma(t)\\
\\
X(t,s)&\mbox{if}&\gamma(t)\leq t \leq t_{k_t+1},\\
\end{array}\right.
\end{equation}
$s \in [t_n,t_{n+1}]$, $t \geq t_{n_0}$ and $n \in [n_0,+\infty
[\cap \mathbb{Z}$.

So,
\[
\begin{array}{rcl}
\psi(t)&=&\displaystyle Z(t,t_{k_t})\left[\Phi(k_t)\Phi(n_0)^{-1}\psi(t_{n_0})\right.\\
\\
&+&\displaystyle\left.\sum_{j=n_0}^{k_t-1} \Phi(k_t)\Phi(j+1)^{-1}\int_{t_j}^{t_{j+1}}\Gamma(t_{j+1},s)f(s)ds\right]\\
\\
&+&\displaystyle\int_{t_{k_t}}^t\Gamma(t,s)f(s)ds.\\
\end{array}
\]
Hence,
\begin{equation}
\label{fvp-depcag}
\begin{array}{rcl}
\psi(t)&=&Z(t,t_{n_0})\psi(t_{n_0})
\\
&+&\displaystyle\sum_{j=n_0}^{k_t-1}\int_{t_{j}}^{t_{j+1}}Z(t,t_{j+1})\Gamma(t_{j+1},s)f(s)ds\\
\\
&+&\displaystyle\int_{t_{k_t}}^t\Gamma(t,s)f(s)ds.\\
\end{array}
\end{equation}

Then, the integral equation (\ref{fvp-depcag}) can be written as
\begin{equation}
\label{fvp-depcag-new}
\psi(t)=Z(t,t_{n_0})\psi(t_{n_0})+\int_{t_0}^t\hat{Z}(t,s)f(s)ds,
\end{equation}
for $t \geq n_0$, where
\begin{equation}
\label{mc} \hat{Z}(t,s)= \left\{\begin{array}{rcl}
Z(t,t_{n+1})\Gamma(t_{n+1},s)&\mbox{if}&t_n \leq s \leq t_{n+1}\\
\\
\Gamma(t,s)&\mbox{if}&t_n\leq s \leq t_{n+1},
\end{array}\right.
\end{equation}

Assume that there is a projection $P:\mathbb{C}^N \to \mathbb{C}^N$ such that
\begin{subequations}
\begin{equation}
\label{26a}
\hat{e} \in (I-P)(\mathbb{C}^N);
\end{equation}
\end{subequations}

Assume that there are a bounded function $h:[t_0,+\infty[ \times
[t_0,+\infty[ \to [0,+\infty[$ and a constant $M>0$ such that
\begin{subequations}
\begin{equation}
\label{27a} \left\|\hat{Z}(t,s)P\right\| \leq
h(t,s)|\tilde{e}(t,s)|,
\end{equation}
for all $t,s \in [t_0,+\infty[$ such that $t \geq s$;
\begin{equation}
\label{27b} \left\|\hat{Z}(t,s)\right\| \leq M|\tilde{e}(t,s)|;
\end{equation}
for all $t,s \in [t_0,+\infty[$ such that $t<s$.

Assume that
\begin{equation}
\label{27c}
h(t,s)\to 0\;\mbox{as}\; t \to +\infty;
\end{equation}
\begin{equation}
\label{27d}
h(t,s) \leq h(t,T)h(T,s)\;\mbox{if}\; t \geq T \geq s.
\end{equation}
\end{subequations}

%Moreover, assume that
%\begin{subequations}
%\begin{equation}
%\label{29a}
%\int_{t_0}^{+\infty} \left\|\tilde{e}(s,t_{k_s})^{-1}R(s)\right\|ds<+\infty.
%\end{equation}
%\end{subequations}

It will be used the variation of parameters formula for DEPCAG  \cite{A04,P2011} as follows.

%Notice that $X(t,s)\hat{e}=\exp\left(\int_s^t
%\hat{\lambda}(\xi)d\xi\right)\hat{e}$ and
%$B(t)\hat{e}=\hat{\lambda}_d(t)\hat{e}$ implies condition
%(\ref{25a}).

%{\tt Formula (\ref{25a}) is naturally obtained under the assumption
%$\xi_n=t_n$ for all $n \in \mathbb{N}_0$.}

Let's extend the conditions (\ref{27a}) and (\ref{27b}) as 
%Then, we have the following dichotomy (\ref{27a}) as
\begin{equation}
\label{27a-nl}
|\hat{Z}(t,s)P| \leq M|\tilde{e}(t,s)|h(t,s),\;\mbox{for}\;t \geq s,
\end{equation}
and
\begin{equation}
\label{27b-nl}
|\hat{Z}(t,s)(I-P)| \leq M|\tilde{e}(t,s)|,\;\mbox{for}\;t \leq s,
\end{equation}
respectively.

Let $\Xi(t,s)=\tilde{e}(t,s)^{-1}\hat{Z}(t,s)$, for $t,s \geq 0$.
Sea $G(t,s): [0,+\infty[^2 \to {\cal M}_{N}(\mathbb{C})$ defined for
\[
G(t,s)=\tilde{e}(t,s)\times\left\{\begin{array}{rcl}
\Xi(t,s)P,&\mbox{if}&t \geq s\\
\\
-\Xi(t,s)(I-P),&\mbox{if}&t<s.
\end{array}
\right.
\]

Then, $\|G(t,s)\| \leq  h(t,s)|\tilde{e}(t,s)|\leq K|\tilde{e}(t,s)|$, for all $t,s \geq 0$.

An additional condition on the perturbation $F(t,x(g(t)))$ is 
\begin{equation}
\label{nl-l1}
\int_{t_0}^{+\infty} |\tilde{e}(t,g(t))|^{-1}\eta(t)dt<+\infty,
\end{equation}

Let $n_0 \in \mathbb{N}_0$. Sea ${\cal B}_{n_0}$ the set of functions $y:[n_0,+\infty[ \to \mathbb{C}^N$ such that $\tilde{e}(\cdot,t_{n_0})^{-1}y \in L^{\infty}$. For $y \in {\cal B}_{n_0}$, sea $\displaystyle \|y\|_{n_0}=\sup_{t \geq t_{n_0}} |\tilde{e}(t,t_{n_0})|^{-1}|y(t)|$. Then,
$({\cal B}_{n_0},\|\cdot\|_{n_0})$ is Banach space, which is isometrically isomorphic to the Banach space $(L^{\infty},\|\cdot\|_{\infty})$.

Let ${\cal N}:{\cal B}_{n_0} \to {\cal B}_{n_0}$ be an operator defined by
\begin{equation}
\label{opabs}
({\cal N}y)(t)=\tilde{e}(t,t_{n_0})\hat{e}+\int_{t_{n_0}}^{+\infty}G(t,s)F(s,g(s))ds,
\end{equation}
for all $y \in {\cal B}_{n_0}$ and $t \geq t_{n_0}$.

Notice that the operator ${\cal N}$ is well defined. In fact, let $y \in {\cal B}_{n_0}$. Then,
\[
\begin{array}{rcl}
|\tilde{e}(t,t_{n_0})|^{-1}|({\cal N}y)(t)| &\leq&\displaystyle |\hat{e}|+M\int_{t_{n_0}}^{t} h(t,s)\left|\tilde{e}(s,g(s))^{-1}F(s,g(s))\right||\tilde{e}(g(s),t_{n_0})^{-1}y(g(s))|ds\\
\\
&+&\displaystyle K\int_{t}^{\infty}\left|\tilde{e}(s,g(s))^{-1}F(s,g(s))\right||\tilde{e}(g(s),t_{n_0})^{-1}y(g(s))|ds\\
\\
 &\leq&
\displaystyle 1+\Theta_{n_0}(t)\|y\|_{n_0},
\end{array}
\]
where $\displaystyle \Theta_{n_0}(t)= M\int_{t_{n_0}}^{t} h(t,s)\left|\tilde{e}(s,g(s))|^{-1}\eta(s)\right|ds+M\int_{t}^{\infty}\left|\tilde{e}(s,g(s))|^{-1}\eta(s)\right|ds$.

By (\ref{nl-l1}), $\displaystyle\Theta_{n_0}(t) \leq \sup_{\tau\geq t_{n_0}}\Theta_{n_0}(\tau)<+\infty$.

So, ${\cal N}y \in {\cal B}_{n_0}$.

Let $n_0$ so large that $\Theta_{n_0}(\infty)<1$. Since
\[
\|{\cal N}y_1-{\cal N}y_2\|_{n_0} \leq \Theta_{n_0}(\infty)\|y_1-y_2\|_{n_0},
\]
for all $y_1, y_2 \in {\cal B}_{n_0}$, by the Banach Fixed Point Theorem, there is a only one $y \in  {\cal B}_{n_0}$ such that $y={\cal N}y$.

Moreover,
\begin{equation}
\label{comparacion}
|\tilde{e}(t,t_{n_0})^{-1} y(t)-\hat{e}| \leq \Theta_{n_0}(t) \|y\|_{n_0}.
\end{equation}

Due to the conditions (\ref{nl-lpchzt}) y (\ref{nl-l1}), we have that
\[
\left|\tilde{e}(\cdot,g(\cdot))|^{-1}\eta(\cdot)\right| \in L^1.
\]

It is not hard to see that $\displaystyle \lim_{t \to +\infty} \Theta_{n_0}(t)=0$. By considering 
\[
w(t)=\tilde{e}(t,t_{n_0})^{-1} y(t)-\hat{e}, 
\]
by the inequality (\ref{comparacion}), we have $w(t) \to 0$ as $t \to +\infty$.

So, the constructed contractive operator allows us to state the following result.
\begin{lema}
\label{nl-lev} Assume for that DEPCAG (\ref{ihomo2}) the conditions (\ref{24a}), (\ref{27c}), (\ref{27d}),
(\ref{27a-nl}), (\ref{27b-nl}) are satisfied and the conditions (\ref{nl-lpchzt}) and (\ref{nl-l1}) are satisfied for $F$.  Then, the
operator ${\cal N}$ defined by (\ref{opabs}) satisfies:
\begin{enumerate}
\item ${\cal N}({\cal B}_{n_0}) \subseteq {\cal B}_{n_0}$;
\item $\tilde{e}(t,n_0)^{-1}\left[({\cal N}y)(t)-\hat{e}\right]\to 0$ as $t \to +\infty$;
\item ${\cal N}$ is contractive for $n_0$ large enough;
\item ${\cal N}$ has a only one fixed point $y_{n_0} \in {\cal B}_{n_0}$, i. e., $y_{n_0}={\cal N}(y_{n_0})$;
\item The fixed point $y_{n_0}$ satisfies the asymptotic formula
\begin{equation}
\label{fas-nl}
y_{n_0}(t)=\tilde{e}(t,n_0)(\hat{e}+w(t)),
\end{equation}
where $w(t) \to 0$ as $t \to +\infty$.
\end{enumerate}
\end{lema}

\section{Main Result}
Now, we are in conditions to present our main result.
\begin{teo}
\label{teoabs} Assume that $\xi_n=t_n$ for all $n \in
\mathbb{N}_0$ an that conditions (\ref{24a}), (\ref{25a}), (\ref{26a}), (\ref{27c}), (\ref{27d}), (\ref{27a-nl}), (\ref{27b-nl}) and (\ref{nl-l1}) are satisfied. Then, the DEPCAG (\ref{idepca2})
 has a solution $y=y(t)$ defined for $t \geq t_{n_0}$ with $n_0$ large enough such that
\begin{equation}
\label{fasintotica-0}
y(t)=\tilde{e}(t,t_{n_0})(\hat{e}+w(t)),
\end{equation}
where $w(t) \to 0$ as $t \to +\infty$.
\end{teo}
\iz{Proof:} Since  the conditions of Teorema \ref{nl-lev} are satisfied, the operator ${\cal N}$ defined by
\[
({\cal N}y)(t)=\tilde{e}(t,t_{n_0})\hat{e}+\int_{t_{n_0}}^{+\infty}G(t,s)F(s,y(g(s)))ds
\]
has a fixed point $y=y(t)$ defined for all $t \geq t_{n_0}$ with $n_0$ large enough with the asymptotic formula (\ref{fas-nl}) which takes the form (\ref{fasintotica-0}).

From (\ref{fvp-depcag-new}) it can be easily proved that $y$
is solution of (\ref{idepca2}). Therefore the  DEPCAG
(\ref{idepca2}) has a solution $y=y(t)$ defined for $t \geq
t_{n_0}$ with $n_0$ large enough with the asymptotic formula (\ref{fasintotica-0})
\begin{flushright}
$\Box$
\end{flushright}
%\iz{Remark:} If we had considered $\xi_n \neq t_n$ as it was assumed, the asymptotic formula (\ref{fasintotica-0}) is
%\begin{equation}
%\label{fasintotica-1}
%y(t)=\tilde{e}(t,t_{n_0})(\hat{e}+w(t)),
%\end{equation}
%where $\tilde{e}$ is given by (\ref{25b}).
%\[
%\begin{array}{rcl}
%E(t,s)&=&e^{\int_s^t \hat{\lambda}(\xi)d\xi}\\
%\\
%&\times&\left(1+e^{\int_{\gamma(t)}^t\hat{\lambda}(\xi)d\xi}\hat{\lambda}_d(u)du\right)\left(1+e^{\int_{\gamma(t)}^{t_{k_t}}\hat{\lambda}(\xi)d\xi}\hat{\lambda}_d(u)du\right)^{-1}\\
%\\
%&\times&\left(\prod_{j=k_s+1}^{k_t}\left[1+\int_{\xi_j}^{t_{j+1}}e^{\int_{\xi_j}^{u}\hat{\lambda}(\xi)d\xi}\hat{\lambda}_d(u)du\right]
%\left[1+\int_{\xi_j}^{t_{j}}e^{\int_{\xi_j}^{u}\hat{\lambda}(\xi)d\xi}\hat{\lambda}_d(u)du\right]^{-1}\right)\\
%\\
%&\times&\left(1+e^{\int_{\gamma(s)}^{t_{k_s+1}}\hat{\lambda}(\xi)d\xi}\hat{\lambda}_d(u)du\right)\left(1+e^{\int_{\gamma(s)}^s\hat{\lambda}(\xi)d\xi}\hat{\lambda}_d(u)du\right)^{-1}.\\
%\end{array}
%\]
%From the proof the Theorem \ref{teoabs}, it can be obtained that DEPCAG (\ref{idepca2}) has a solution with asymptotic formula (\ref{fasintotica-1}). In similar way, the DEPCAG (\ref{idepca2}) can be changed by the DEPCAG:
%\begin{equation}
%\label{depcag3}
%\frac{dy}{dt}=A(t)y(t)+B(t)y(\gamma(t))+F(t,g(t)),
%\end{equation}
%where $F$ satisfies (\ref{nl-lpchzt}) and (\ref{nl-l1}).  Then the DEPCAG (\ref{depcag3}) has a solution with asymptotic formula (\ref{fasintotica-1}).

\iz{Example}
Consider the diagonal matrices in ${\cal M}_N(\mathbb{C})$,
\[
\Lambda_A(t)=\mbox{diag}(a_1(t),a_2(t),\ldots,a_N(t))\;\mbox{and}\; \Lambda_B(t)=\mbox{diag}(b_1(t),b_2(t),\ldots,b_N(t))
\]
DEPCAG with delayed piecewise constant argument $\xi_n=t_n$, for all $n \in \mathbb{N}_0$,
\begin{equation}
\label{depcag1}
z'(t)=\Lambda_A(t)z(t)+\Lambda_B(t)z(\gamma(t)),
\end{equation}
for $t \geq t_0$.
Then, the solutions of (\ref{depcag1}) can be written as
\[
y(t)=Z(t,s)y(s),
\]
where
$\displaystyle Z(t,s)=\mbox{diag}(e_1(t,s),e_2(t,s),\ldots,e_N(t,s))$,
\[
\begin{array}{rcl}
e_l(t,s)&=&e^{\int_s^ta_l(\xi)d\xi}\left(1+\int_{t_{\gamma(s)}}^s e^{-\int_{\gamma(s)}^{\sigma} a_l(\xi)d\xi}b_l(\sigma)d\sigma\right)^{-1}\\
&\times&\left[\prod_{m=\gamma(s)}^{k_t-1} \left(1+\int_{t_m}^{t_{m+1}} e^{-\int_{\gamma(s)}^{\sigma} a_l(\xi)d\xi}b_l(\sigma)d\sigma \right)\right]\\
&\times&\left(1+\int_{\gamma(t)}^t e^{-\int_{\gamma(t)}^{\sigma} a_l(\xi)d\xi}b_l(\sigma)d\sigma\right),
\end{array}
\]
for $l \in \{1,\ldots,N\}$ y $t \geq s$.

We assume that
\begin{equation}
\label{hinv}
1+\int_{t_m}^{t} e^{-\int_{\gamma(s)}^{\sigma} a_l(\xi)d\xi}b_l(\sigma)d\sigma \neq 0,\;\mbox{for all}\;l \in \{1,\ldots,N\},\;t \in [t_m,t_{m+1}[\;\mbox{y}\; m \in \mathbb{N}_0.
\end{equation}
This is equivalent to (\ref{24a}) an guarantees that $Z(t,s)$ is invertible for $t \geq s$ and the existence of $Z(t,s)$ for $t<s$. Moreover, $e_l(t,s)=\frac{1}{e_l(s,t)}$ for $t<s$.

The following is a direct application of theorem \ref{teoabs}.
\begin{co} Assume that there is $k \in \{1,\ldots,N\}$ such that
\begin{enumerate}[(a)]
\item $\displaystyle \lim_{t \to +\infty} \left|\frac{e_l(t,s)}{e_k(t,s)}\right|=0$, for $l<k$,
\item There is $C \geq 0$ such that $\left|\frac{e_l(t,s)}{e_k(t,s)}\right| \leq C$, for $l \geq k$.
\end{enumerate}
Let $\{R(t)\}_{t \geq 0}$ a family in ${\cal M}_N(\mathbb{C})$ such that $R(\cdot)$ is locally integrable and
\begin{equation}
\label{l1-01}
\sum_{n=n_0}^{+\infty} \int_{t_n}^{t_{n+1}} |e_k(s,t_n)^{-1}|\|R(s)\|ds\;\mbox{is convergent.}
\end{equation}
Then, the DEPCAG
\[
y'(t)=\Lambda_A(t)y(t)+\Lambda_B(t)y(\gamma(t))+R(t)y(\gamma(t)),
\]
has a solution $y=y(t)$ defined for $t \geq t_{n_0}$ with $n_0$ large enough such that
\begin{equation}
\label{fasintotica}
y(t)=e_k(t,t_{n_0})(e_k+w(t)),
\end{equation}
where $e_k$ is the $k$-th vector of the canonical base of $\mathbb{C}^N$ and $w(t) \to 0$ as $t \to +\infty$.
\end{co}

\end{document}